%% file: LinksDeWit2000.tex
\documentclass[fleqn,twoside,10pt,a4paper]{article}

\usepackage{amssymb}
\usepackage{epsfig}

\title{
  Link Invariants Associated with Gauge Equivalent Solutions of the
  Yang--Baxter Equation:  the One-Parameter Family of Minimal Typical
  Representations of {$U_q[gl(2|1)]$}
}

\author{
  Jon R Links%
  \footnote{
    Centre for Mathematical Physics, Department of Mathematics,
    The University of Queensland, 4072,
    Australia.
    \texttt{jrl@maths.uq.edu.au}
  }
  \\
  and
  \\
  David~~De Wit%
  \footnote{
    RIMS, Kyoto University 606-8502, Japan.
    \texttt{ddw@kurims.kyoto-u.ac.jp}
  }
}

\begin{document}

\maketitle

\begin{abstract}
  In this paper we investigate the construction of state models for
  link invariants using representations of the braid group obtained
  from various \emph{gauge} choices for a solution of the trigonometric
  Yang--Baxter equation. Our results show that it is possible to obtain
  invariants of \emph{regular isotopy} (as defined by Kauffman) which
  may not be \emph{ambient isotopic}. We illustrate our results with
  explicit computations using solutions of the trigonometric
  Yang--Baxter equation associated with the one-parameter family of
  minimal typical representations of the quantum superalgebra
  $U_q[gl(2|1)]$. We have implemented \textsc{Mathematica} code to
  evaluate the invariants for all prime knots up to $10$ crossings.
\end{abstract}

%%%%%%%%%%%%%%%%%%%%%%%%%%%%%%%%%%%%%%%%%%%%%%%%%%%%%%%%%%%%%%%%%%%%%%%%

\section{Introduction}

The close connection between solutions of the quantum Yang--Baxter
equation (QYBE) and the evaluation of invariants for oriented knots and
links through representations of the braid group is well known
\cite{YangGe:89}.  The algebraic properties of quantum algebras and
superalgebras provide systematic means to construct solutions of the
QYBE which can in turn be used to explicitly compute these invariants.
The seminal example of this method is the use of the six vertex
solution of the QYBE in the evaluation of the Jones polynomial invariant
\cite{Kauffman:93}.

The spectral parameter dependent solutions of the trigonometric
Yang--Baxter equation (TYBE) arise as evaluation (or loop)
representations of \emph{affine} quantum algebras and
superalgebras.  Braid group representations are obtained from these
solutions by taking the limit as the spectral parameter approaches
infinity.  We refer to these limiting cases as
\emph{quantum R matrices} which satisfy the QYBE. As demonstrated by
Bracken et al.  \cite{BrackenGouldZhangDelius:94a}, gauge equivalent
solutions of the TYBE are obtained by considering different gradations
of the underlying affine (super)algebraic structure. The significant
point is that the explicit braid group representation obtained in the
infinite spectral parameter limit depends on the choice of the
gradation. This feature has already been observed in
\cite{BrackenGouldZhangDelius:94a}. This means that given a single
solution of the TYBE, there are many possibilities to in turn obtain a
solution of the QYBE.

Our aim in this work is to investigate the construction of link
invariants for the case of solutions of the QYBE obtained through
non-standard choices of gradation (i.e., not the \emph{homogeneous}
gradation) for the affine (super)algebraic structure which underlies
the TYBE solutions. We find that in general we are only able to define
invariants of \emph{regular isotopy}, which is to say that the results
are invariant under the second and third Reidemeister moves. However,
we will show that in our examples we can relate the regular isotopy
invariants to well known invariants of \emph{ambient isotopy}, viz
invariant under all three Reidemeister moves. (Our terminology for
regular and ambient isotopy is adopted from \cite{Kauffman:93}.)

Our specific calculations are performed using the solution of the TYBE
associated with the one-parameter family of minimal typical
representations of the quantum superalgebra $U_q[gl(2|1)]$. For the
choice of the homogeneous gradation for the untwisted affine extension
$U_q[gl(2|1)^{(1)}]$, one obtains the Links--Gould invariants recently
investigated in detail in
\cite{DeWit:99a,DeWitKauffmanLinks:99a}.
However, for different choices of the gradation the process yields
other invariants which can be related back to both the Jones and
Alexander--Conway polynomial invariants.

We begin the paper with a description of gauge equivalent solutions of
the TYBE without appealing to the gradation structure of quantum affine
(super)algebras. For the sake of simplicity, we present these solutions
in terms of a simple basis transformation satisfying some particular
constraints.

%%%%%%%%%%%%%%%%%%%%%%%%%%%%%%%%%%%%%%%%%%%%%%%%%%%%%%%%%%%%%%%%%%%%%%%%

\section{Gauge equivalent solutions of the TYBE}

For an arbitrary vector space $V$, let
$R(x) \in \mathrm{End} (V\otimes V)$ satisfy the TYBE on the
tensor product space $V \otimes V \otimes V$:
\begin{eqnarray}
  R_{12}(x)
  R_{13}(xy)
  R_{23}(y)
  =
  R_{23}(y)
  R_{13}(xy)
  R_{12}(x).
  \label{eq:TYBE}
\end{eqnarray}
where $R_{12}(x) \triangleq R(x) \otimes I$, etc, and
$x$ and $y$ are arbitrary complex parameters.  Consider an
invertible matrix $A(x) \in \mathrm{End} (V)$, with the properties:
\begin{eqnarray}
  \left.
  \begin{array}{rcl}
    A(x) A(y)
    & = &
    A(x y)
    \\
    \left[
      R(x), \, A_1(y) A_2(y)
    \right]
    & = &
    0,
  \end{array}
  \qquad \qquad \qquad \qquad \qquad \qquad \qquad \qquad
  \right\}
  \label{eq:AProperties}    
\end{eqnarray}
where $A_1(x) \triangleq A(x) \otimes I $ and
$ A_2(x) \triangleq I \otimes A(x)$.
We immediately deduce that:
\begin{eqnarray*}
  A(1)
  & = &
  I,
  \\
  A(x)A(y)
  & = &
  A(y)A(x),
  \\
  A^{-1}(x)
  & = &
  A(\overline{x}),
\end{eqnarray*}
where throughout the paper, we liberally write
$\overline{X}$ for $X^{-1}$, for various $X$.

Now set $\mathcal{R}(x)=A_1(x)R(x)A_1(\overline{x}).$
It is an algebra exercise to prove:
\begin{eqnarray*}
  \begin{array}{rcl}
    & &
    \hspace{-20pt}
    \mathcal{R}_{12}(x)\mathcal{R}_{13}(xy)\mathcal{R}_{23}(y)
    \\
    & & =
    A_1(x)R_{12}(x)A_1(\overline{x})A_1(xy)
    R_{13}(xy)A_1(\overline{x}\overline{y})A_2(y)
    R_{23}(y)A_2(\overline{y})
    \\
    & & =
    A_1(x)R_{12}(x)A_1(y)A_2(y)R_{13}(xy)
    R_{23}(y)A_1(\overline{x}\overline{y})A_2(\overline{y}) \\
    & & =
    A_1(xy)A_2(y)R_{12}(x)R_{13}(xy)
    R_{23}(y)A_1(\overline{x}\overline{y})A_2(\overline{y})  \\
    & & =
    A_1(xy)A_2(y)R_{23}(y)R_{13}(xy)
    R_{12}(x)A_1(\overline{x}\overline{y})A_2(\overline{y})  \\
    & & =
    A_2(y)R_{23}(y)A_1(xy)
    R_{13}(xy)A_1(\overline{y})A_2(\overline{y})
    R_{12}(x)A_1(\overline{x}) 
    \\
    & & =
    A_2(y)R_{23}(y)A_2(\overline{y})A_1(xy)
    R_{13}(xy)A_1(\overline{y})R_{12}(x)A_1(\overline{x}) 
    \\
    & & =
    A_2(y)R_{23}(y)A_2(\overline{y})A_1(xy)
    R_{13}(xy)A_1(\overline{x}\overline{y})A_1(x)
    R_{12}(x)A_1(\overline{x}) \\
    & & =
    \mathcal{R}_{23}(y)\mathcal{R}_{13}(xy)\mathcal{R}_{12}(x).
  \end{array}
\end{eqnarray*}

The above calculation shows that $\mathcal{R}(x)$ also satisfies
(\ref{eq:TYBE}).  The connection between $R(x)$ and $\mathcal{R}(x)$ in
terms of the gradation of the underlying affine algebraic structure for
these solutions is discussed in \cite{BrackenGouldZhangDelius:94a}.
For our purposes, we need not describe this in detail, and simply refer
to $R(x)$ and $\mathcal{R}(x)$ as \emph{gauge equivalent}.
The essential point is that the limit of $\mathcal{R}(x)$ as $x\to
\infty$ yields different quantum R matrices depending on the choice of
$A(x)$.

%%%%%%%%%%%%%%%%%%%%%%%%%%%%%%%%%%%%%%%%%%%%%%%%%%%%%%%%%%%%%%%%%%%%%%%%

\section{Trigonometric R matrix with gauge parameters}

Applying the above to a particular example, the following trigonometric
R matrix $\check{\mathcal{R}}^{r,s}(u)$ with \emph{gauge parameters}
$r$ and $s$ arises from the one-parameter family of minimal typical
representations of $U_q[gl(2|1)]$.  (Here, we have replaced variable
$x$ with $u$, defined by $x\triangleq q^u$.) The operator
$\check{\mathcal{R}}^{r,s}(u)$ has $36$ nonzero components, and is
scaled such that its first component is $1$.  It satisfies the TYBE in
the additive form:
\begin{eqnarray}
  \check{\mathcal{R}}_{12}(u)
  \check{\mathcal{R}}_{23}(u+v)
  \check{\mathcal{R}}_{12}(v)
  =
  \check{\mathcal{R}}_{23}(v)
  \check{\mathcal{R}}_{12}(u+v)
  \check{\mathcal{R}}_{23}(u).
  \label{eq:braidTYBE}
\end{eqnarray} 
(Confirming this involves manipulating expressions in a total of $6$
variables, viz: the representation variables $q$ and $\alpha$, the
spectral variables $u$ and $v$, and the gauge parameters $r$ and $s$.)
Explicitly, $\check{\mathcal{R}}^{r,s}(u)$ is:
\begin{eqnarray*}
  &&
  \hspace{-45pt}
  \left\{
  \begin{array}{@{\hspace{0mm}}c@{\hspace{0mm}}}
    e^{1 1}_{1 1}
  \end{array}
  \right\},
  \qquad
  \frac{
    [\alpha + u]
  }{
    [\alpha - u]
  }
  \left\{
  \begin{array}{@{\hspace{0mm}}c@{\hspace{0mm}}}
    e^{2 2}_{2 2}
    \\
    e^{3 3}_{3 3}
  \end{array}
  \right\},
  \qquad
  \frac{
    [\alpha + u]
    [1 + \alpha + u]
  }{
    [\alpha - u]
    [1 + \alpha - u]
  }
  \left\{
  \begin{array}{@{\hspace{0mm}}c@{\hspace{0mm}}}
    e^{4 4}_{4 4}
  \end{array}
  \right\},
  \\
  &&
  \hspace{-45pt}
  \frac{
    [\alpha]
  }{
    [\alpha - u]
  }
  \left\{
  \begin{array}{@{\hspace{0mm}}c@{\hspace{0mm}}}
    r^u
    \overline{q}^u
    e^{1 2}_{1 2},
    s^u
    \overline{q}^u
    e^{1 3}_{1 3}
    \\
    \overline{r}^u
    q^u
    e^{2 1}_{2 1},
    \overline{s}^u
    q^u
    e^{3 1}_{3 1}
  \end{array}
  \right\},
  \qquad
  \frac{
    [\alpha]
    [1 + \alpha]
  }{
    [\alpha - u]
    [1 + \alpha - u]
  }
  \left\{
  \begin{array}{@{\hspace{0mm}}c@{\hspace{0mm}}}
    r^u
    s^u
    \overline{q}^{2 u}
    e^{1 4}_{1 4}
    \\
    \overline{r}^u
    \overline{s}^u
    q^{2 u}
    e^{4 1}_{4 1}
  \end{array}
  \right\},
  \\
  &&
  \hspace{-45pt}
  \frac{
    1
  }{
    \Delta^2
    [\alpha - u]
    [1 + \alpha - u]
  }
  \left\{
    \begin{array}{@{\hspace{0mm}}c@{\hspace{0mm}}}
      \overline{r}^u
      s^u
      f(\overline{q})
      e^{2 3}_{2 3}
      \\
      r^u
      \overline{s}^u
      f(q)
      e^{3 2}_{3 2}
    \end{array}
  \right\},
  \qquad
  \frac{
    [1 + \alpha]
    [\alpha + u]
  }
    {{[\alpha - u]}
    [1 + \alpha - u]
  }
  \left\{
    \begin{array}{@{\hspace{0mm}}c@{\hspace{0mm}}}
      s^u
      \overline{q}^u
      e^{2 4}_{2 4},
      r^u
      \overline{q}^u
      e^{3 4}_{3 4}
      \\
      \overline{s}^u
      q^u
      e^{4 2}_{4 2},
      \overline{r}^u
      q^u
      e^{4 3}_{4 3}
    \end{array}
  \right\},
  \\
  &&
  \hspace{-45pt}
  -
  \frac{
    [u]
  }{
    [\alpha - u]
  }
  \left\{
  \begin{array}{@{\hspace{0mm}}c@{\hspace{0mm}}}
    e^{1 2}_{2 1},
    e^{1 3}_{3 1}
    \\
    e^{2 1}_{1 2},
    e^{3 1}_{1 3}
  \end{array}
  \right\},
  \qquad
  -
  \frac{
    [1 - u]
    [u]
  }{
    [\alpha - u]
    [1 + \alpha - u]
  }
  \left\{
  \begin{array}{@{\hspace{0mm}}c@{\hspace{0mm}}}
    e^{1 4}_{4 1}
    \\
    e^{4 1}_{1 4}
  \end{array}
  \right\},
  \\
  &&
  \hspace{-45pt}
  -
  \frac{
    [u]^2
  }{
    [\alpha - u]
    [1 + \alpha - u]
  }
  \left\{
  \begin{array}{@{\hspace{0mm}}c@{\hspace{0mm}}}
    e^{2 3}_{3 2}
    \\
    e^{3 2}_{2 3}
  \end{array}
  \right\},
  \qquad
  \frac{
    [u]
    [\alpha + u]
  }{
    [\alpha - u]
    [1 + \alpha - u]
  }
  \left\{
  \begin{array}{@{\hspace{0mm}}c@{\hspace{0mm}}}
    e^{2 4}_{4 2},
    e^{3 4}_{4 3}
    \\
    e^{4 2}_{2 4},
    e^{4 3}_{3 4}
  \end{array}
  \right\},
  \\
  &&
  \hspace{-48pt}
  \frac{
    {{[\alpha]}}^{\frac{1}{2}}
    {{[1 + \alpha]}}^{\frac{1}{2}}
    [u]
  }{
    [\alpha - u]
    [1 + \alpha - u]
  }
  \left\{
  \begin{array}{@{\hspace{0mm}}c@{\hspace{0mm}}}
    +
    {\displaystyle \frac{r^u}{q^u}}
    q^{\frac{1}{2}}
    \left\{
    \begin{array}{@{\hspace{0mm}}c@{\hspace{0mm}}}
      e^{1 4}_{3 2}
      \\
      e^{3 2}_{1 4}
    \end{array}
    \right\},
    -
    {\displaystyle \frac{q^u}{r^u}}
    \overline{q}^{\frac{1}{2}}
    \left\{
    \begin{array}{@{\hspace{0mm}}c@{\hspace{0mm}}}
      e^{4 1}_{2 3}
      \\
      e^{2 3}_{4 1}
    \end{array}
    \right\},
    +
    {\displaystyle \frac{q^u}{s^u}}
    q^{\frac{1}{2}}
    \left\{
    \begin{array}{@{\hspace{0mm}}c@{\hspace{0mm}}}
      e^{3 2}_{4 1}
      \\
      e^{4 1}_{3 2}
    \end{array}
    \right\},
    -
    {\displaystyle \frac{s^u}{q^u}}
    \overline{q}^{\frac{1}{2}}
    \left\{
    \begin{array}{@{\hspace{0mm}}c@{\hspace{0mm}}}
      e^{2 3}_{1 4}
      \\
      e^{1 4}_{2 3}
    \end{array}
    \right\}
  \end{array}
  \right\},
\end{eqnarray*}
where
$
  f(q)
  \triangleq
  - 2 q
  + q^{2 u}
  (q - \overline{q})
  +              q^{1 + 2 \alpha}
  +   \overline{q}^{1 + 2 \alpha}
$
and we have defined $ \Delta \triangleq q - \overline{q}$, as well as
writing $[x]$ for $(q^x - \overline{q}^x)/(q-\overline{q})$.

This solution of the TYBE originates in the following trigonometric R
matrix $\check{R}(u)\equiv \check{\mathcal{R}}^{r=1,s=1}(u)$, which may
be regarded as a \emph{gauge-free} version of
$\check{\mathcal{R}}^{r,s}(u)$.
\begin{eqnarray*}
  &&
  \hspace{-45pt}
  \left\{
  \begin{array}{@{\hspace{0mm}}c@{\hspace{0mm}}}
    e^{1 1}_{1 1}
  \end{array}
  \right\},
  \qquad
  \frac{
    [\alpha + u]
  }{
    [\alpha - u]
  }
  \left\{
  \begin{array}{@{\hspace{0mm}}c@{\hspace{0mm}}}
    e^{2 2}_{2 2}
    \\
    e^{3 3}_{3 3}
  \end{array}
  \right\},
  \qquad
  \frac{
    [\alpha + u]
    [1 + \alpha + u]
  }{
    [\alpha - u]
    [1 + \alpha - u]
  }
  \left\{
  \begin{array}{@{\hspace{0mm}}c@{\hspace{0mm}}}
    e^{4 4}_{4 4}
  \end{array}
  \right\},
  \\
  &&
  \hspace{-45pt}
  \frac{
    [\alpha]
  }{
    [\alpha - u]
  }
  \left\{
  \begin{array}{@{\hspace{0mm}}c@{\hspace{0mm}}}
    \overline{q}^u
    e^{1 2}_{1 2},
    \overline{q}^u
    e^{1 3}_{1 3}
    \\
    q^u
    e^{2 1}_{2 1},
    q^u
    e^{3 1}_{3 1}
  \end{array}
  \right\},
  \qquad
  \frac{
    [\alpha]
    [1 + \alpha]
  }{
    [\alpha - u]
    [1 + \alpha - u]
  }
  \left\{
  \begin{array}{@{\hspace{0mm}}c@{\hspace{0mm}}}
    \overline{q}^{2 u}
    e^{1 4}_{1 4}
    \\
    q^{2 u}
    e^{4 1}_{4 1}
  \end{array}
  \right\},
  \\
  &&
  \hspace{-45pt}
  \frac{
    1
  }{
    \Delta^2
    [\alpha - u]
    [1 + \alpha - u]
  }
  \left\{
    \begin{array}{@{\hspace{0mm}}c@{\hspace{0mm}}}
      f(\overline{q})
      e^{2 3}_{2 3}
      \\
      f(q)
      e^{3 2}_{3 2}
    \end{array}
  \right\},
  \qquad
  \frac{
    [1 + \alpha]
    [\alpha + u]
  }
    {{[\alpha - u]}
    [1 + \alpha - u]
  }
  \left\{
    \begin{array}{@{\hspace{0mm}}c@{\hspace{0mm}}}
      \overline{q}^u
      e^{2 4}_{2 4},
      \overline{q}^u
      e^{3 4}_{3 4}
      \\
      q^u
      e^{4 2}_{4 2},
      q^u
      e^{4 3}_{4 3}
    \end{array}
  \right\},
  \\
  &&
  \hspace{-45pt}
  -
  \frac{
    [u]
  }{
    [\alpha - u]
  }
  \left\{
  \begin{array}{@{\hspace{0mm}}c@{\hspace{0mm}}}
    e^{1 2}_{2 1},
    e^{1 3}_{3 1}
    \\
    e^{2 1}_{1 2},
    e^{3 1}_{1 3}
  \end{array}
  \right\},
  \qquad
  -
  \frac{
    [1 - u]
    [u]
  }{
    [\alpha - u]
    [1 + \alpha - u]
  }
  \left\{
  \begin{array}{@{\hspace{0mm}}c@{\hspace{0mm}}}
    e^{1 4}_{4 1}
    \\
    e^{4 1}_{1 4}
  \end{array}
  \right\},
  \\
  &&
  \hspace{-45pt}
  -
  \frac{
    [u]^2
  }{
    [\alpha - u]
    [1 + \alpha - u]
  }
  \left\{
  \begin{array}{@{\hspace{0mm}}c@{\hspace{0mm}}}
    e^{2 3}_{3 2}
    \\
    e^{3 2}_{2 3}
  \end{array}
  \right\},
  \qquad
  \frac{
    [u]
    [\alpha + u]
  }{
    [\alpha - u]
    [1 + \alpha - u]
  }
  \left\{
  \begin{array}{@{\hspace{0mm}}c@{\hspace{0mm}}}
    e^{2 4}_{4 2},
    e^{3 4}_{4 3}
    \\
    e^{4 2}_{2 4},
    e^{4 3}_{3 4}
  \end{array}
  \right\},
  \\
  &&
  \hspace{-45pt}
  \frac{
    {{[\alpha]}}^{\frac{1}{2}}
    {{[1 + \alpha]}}^{\frac{1}{2}}
    [u]
  }{
    [\alpha - u]
    [1 + \alpha - u]
  }
  \left\{
  \begin{array}{@{\hspace{0mm}}c@{\hspace{0mm}}}
    +
    \overline{q}^{u - \frac{1}{2}}
    \left\{
    \begin{array}{@{\hspace{0mm}}c@{\hspace{0mm}}}
      e^{1 4}_{3 2}
      \\
      e^{3 2}_{1 4}
    \end{array}
    \right\},
    -
    q^{u - \frac{1}{2}}
    \left\{
    \begin{array}{@{\hspace{0mm}}c@{\hspace{0mm}}}
      e^{4 1}_{2 3}
      \\
      e^{2 3}_{4 1}
    \end{array}
    \right\},
    +
    q^{u + \frac{1}{2}}
    \left\{
    \begin{array}{@{\hspace{0mm}}c@{\hspace{0mm}}}
      e^{3 2}_{4 1}
      \\
      e^{4 1}_{3 2}
    \end{array}
    \right\},
    -
    \overline{q}^{u + \frac{1}{2}}
    \left\{
    \begin{array}{@{\hspace{0mm}}c@{\hspace{0mm}}}
      e^{2 3}_{1 4}
      \\
      e^{1 4}_{2 3}
    \end{array}
    \right\}
  \end{array}
  \right\}.
\end{eqnarray*}

The above solution of (\ref{eq:braidTYBE}) 
is obtained from the representation 
theory of the quantum superalgebra $U_q[gl(2|1)]$ and the \emph{Baxterization}
procedure (e.g. see \cite{DeliusGouldLinksZhang:95a}).
It is related to 
$R(u)$ through a transposition:
\begin{eqnarray*}
  \check{R}(u)^{a c}_{b d}
  =
  R(u)^{c a}_{b d}.
\end{eqnarray*}
From it we may construct a quantum R matrix
$\check{R}=\lim_{u\to\infty}\check{R}(u)$, and build a link invariant
\cite{DeWit:99a,DeWitKauffmanLinks:99a}.

By the following procedure, $\check{\mathcal{R}}^{r,s}(u)$ may be
obtained from $\check{R}(u)$.  Let us define:
\begin{eqnarray}
  {\mathcal{R}}(u)^{a c}_{b d}
  \triangleq
  A(u)^{a}_{e}
  \cdot
  R(u)^{e c}_{f d}
  \cdot
  A(-u)^{f}_{b},
  \label{eq:DefnofcalTrigR}
\end{eqnarray}
where $A(u)$ is a $4\times 4$ matrix satisfying the properties
(\ref{eq:AProperties}) which ensure that $\mathcal{R}(u)$ also
satisfies (\ref{eq:TYBE}).  Specifically, we will use the following
diagonal $A(u)\equiv A^{r,s}(u)$ containing gauge parameters $r$ and
$s$:
\begin{eqnarray}
  A^{r,s}(u)
  \triangleq
  \mathrm{diag}
  \left\{
    1, r^{u}, s^{u}, r^{u} s^{u}
  \right\},
  \label{eq:Aisdiagonal}
\end{eqnarray}
thus we will write $\mathcal{R}^{r,s}(u)$ rather than
$\mathcal{R}(u)$.  We define $\check{\mathcal{R}}^{r,s}(u)$ by again
transposing:
\begin{eqnarray}
  \check{\mathcal{R}}^{r,s}(u)^{a c}_{b d}
  \triangleq
  \mathcal{R}^{r,s}(u)^{c a}_{b d}.
  \label{eq:DefnofcheckcalRrs}
\end{eqnarray}
This $\check{\mathcal{R}}^{r,s}(u)$ satisfies (\ref{eq:braidTYBE}), and
indeed is the object which we first introduced.

To investigate the possibility of constructing link invariants, we
will select appropriate gauge choices
$r \equiv r_i(q)$, $s \equiv s_i(q)$
and take the spectral limit $u\to\infty$ of
$\check{\mathcal{R}}^{r_i,s_i}(u)$ to yield a quantum R matrix 
$\check{R}^i$:
\begin{eqnarray}
  \check{R}^{i}
  \triangleq
  \lim_{u\to\infty}
    \check{\mathcal{R}}^{r_i,s_i}(u),
  \label{eq:QRMxDefn}
\end{eqnarray}
which satisfies the (non-parametric) QYBE:
\begin{eqnarray}
  \check{R}_{12}
  \check{R}_{23}
  \check{R}_{12}
  =
  \check{R}_{23}
  \check{R}_{12}
  \check{R}_{23},
  \label{eq:QYBEbraid}
\end{eqnarray}
which is graphically depictable in terms of braids.
(In the above, $i$ is an index to number the gauge choice, and the
functions $r_i$ and $s_i$ may also contain variable $q$.)

Substituting (\ref{eq:DefnofcalTrigR}) and (\ref{eq:Aisdiagonal}) into
(\ref{eq:DefnofcheckcalRrs}) allows us to write (no sums on $b$ and
$c$):
\begin{eqnarray*}
  \check{R}^{r,s}(u)^{a c}_{b d}
  =
  A^{r,s}(u)^{c}_{c}
  \cdot
  \check{R}(u)^{a c}_{b d}
  \cdot
  A^{r,s}(-u)^{b}_{b},
\end{eqnarray*}
hence we have for grading choice $i$ a quantum R matrix $R^i$ with
components:
\begin{eqnarray*}
  (\check{R}^i)^{a c}_{b d}
  =
  \lim_{u\to\infty}
  A^{r_i,s_i}(u)^{c}_{c}
  \cdot
  \check{R}(u)^{a c}_{b d}
  \cdot
  A^{r_i,s_i}(-u)^{b}_{b}.
\end{eqnarray*}
The gauge choices that we shall make are depicted in Table
\ref{tab:gauges}.

\begin{table}[ht]
  \centering
  \begin{tabular}{ccc}
    Case ($i$) & $r_i,s_i$ & $A^{r^i,s^i}(u)$
    \\[1mm]
    \cline{1-1} \cline{2-2} \cline{3-3}
    \\[-3mm]
    1 & $r_1=1$, $s_1=1$ & $ I $,
    \\
    2 & $r_2=1$, $s_2=q$ &
      $\mathrm{diag} \left\{ 1, 1, q^{u}, q^{u} \right\}$
    \\
    3 & $r_3=q$, $s_3=q$ &
      $\mathrm{diag} \left\{ 1, q^{u}, q^{u}, q^{2u} \right\}$
    \\
    4 & $r_4 s_4=q^{2}$, $s_4>r_4>1$ &
      $
        \mathrm{diag}
        \left\{ 1, r^{u}, r^{u} \overline{q}^u, q^{2u} \right\}
      $
  \end{tabular}
  \caption{Gauge choices.}
  \label{tab:gauges}
\end{table}
From Case 1, we recover the ungauged situation.

%%%%%%%%%%%%%%%%%%%%%%%%%%%%%%%%%%%%%%%%%%%%%%%%%%%%%%%%%%%%%%%%%%%%%%%%

\section{Quantum R matrices}

In the spectral limit $u\to \infty$, our trigonometric R matrix
$\check{\mathcal{R}}^{r_i,s_i}(u)$ becomes a quantum R matrix
$\check{R}^i$ in variables $q$ and $\alpha$, where $i$ is the index of
the gauge choice (\ref{eq:QRMxDefn}).  Here, we present the
$\check{R}^i$, in terms of `internal' variables, $p\triangleq
q^{\alpha+\frac{1}{2}}$ and $Q\triangleq q^{\frac{1}{2}}$, which
simplify computations.

$\check{R}^1$ has $26$ nonzero components:
\begin{eqnarray*}
  &&
  \left\{
  \begin{array}{@{\hspace{0mm}}c@{\hspace{0mm}}}
    e^{1 1}_{1 1}
  \end{array}
  \right\},
  \qquad
  -
  p^2 \overline{Q}^2
  \left\{
  \begin{array}{@{\hspace{0mm}}c@{\hspace{0mm}}}
    e^{2 2}_{2 2}
    \\
    e^{3 3}_{3 3}
  \end{array}
  \right\},
  \qquad
  p^4
  \left\{
  \begin{array}{@{\hspace{0mm}}c@{\hspace{0mm}}}
    e^{4 4}_{4 4}
  \end{array}
  \right\},
  \\
  &&
  -
  p \overline{Q}
  (p \overline{Q} - \overline{p} Q)
  \left\{
  \begin{array}{@{\hspace{0mm}}c@{\hspace{0mm}}}
    e^{2 1}_{2 1}
    \\
    e^{3 1}_{3 1}
  \end{array}
  \right\},
  \qquad
  p^3 \overline{Q}
  (p Q - \overline{p} \overline{Q})
  \left\{
  \begin{array}{@{\hspace{0mm}}c@{\hspace{0mm}}}
    e^{4 2}_{4 2}
    \\
    e^{4 3}_{4 3}
  \end{array}
  \right\},
  \\
  &&
  p^2
  (p \overline{Q} - \overline{p} Q)
  (p Q - \overline{p} \overline{Q})
  \left\{
  \begin{array}{@{\hspace{0mm}}c@{\hspace{0mm}}}
    e^{4 1}_{4 1}
  \end{array}
  \right\},
  \qquad
  p^2 ( Q^{2} - \overline{Q}^{2} )
  \left\{
  \begin{array}{@{\hspace{0mm}}c@{\hspace{0mm}}}
    e^{3 2}_{3 2}
  \end{array}
  \right\},
  \\
  &&
  p \overline{Q}
  \left\{
  \begin{array}{@{\hspace{0mm}}c@{\hspace{0mm}}}
    e^{1 2}_{2 1},
    e^{1 3}_{3 1}
    \\
    e^{2 1}_{1 2},
    e^{3 1}_{1 3}
  \end{array}
  \right\},
  \qquad
  p^3 \overline{Q}
  \left\{
  \begin{array}{@{\hspace{0mm}}c@{\hspace{0mm}}}
    e^{2 4}_{4 2},
    e^{3 4}_{4 3}
    \\
    e^{4 2}_{2 4},
    e^{4 3}_{3 4}
  \end{array}
  \right\},
  \qquad
  p^2 \overline{Q}^2
  \left\{
  \begin{array}{@{\hspace{0mm}}c@{\hspace{0mm}}}
    e^{1 4}_{4 1}
    \\
    e^{4 1}_{1 4}
  \end{array}
  \right\},
  \qquad
  -
  p^2
  \left\{
  \begin{array}{@{\hspace{0mm}}c@{\hspace{0mm}}}
    e^{2 3}_{3 2}
    \\
    e^{3 2}_{2 3}
  \end{array}
  \right\},
  \\
  &&
  p^2
  {(p\overline{Q} - \overline{p}Q)}^{\frac{1}{2}}
  {(p Q - \overline{p} \overline{Q})}^{\frac{1}{2}}
  \left\{
  \begin{array}{@{\hspace{0mm}}c@{\hspace{0mm}}}
    -\overline{Q}
    \left\{
    \begin{array}{@{\hspace{0mm}}c@{\hspace{0mm}}}
      e^{4 1}_{2 3}
      \\
      e^{2 3}_{4 1}
    \end{array}
    \right\},
    Q
    \left\{
    \begin{array}{@{\hspace{0mm}}c@{\hspace{0mm}}}
      e^{3 2}_{4 1}
      \\
      e^{4 1}_{3 2}
    \end{array}
    \right\}
  \end{array}
  \right\}.
\end{eqnarray*}

$\check{R}^2$ has $20$ nonzero components:
\begin{eqnarray*}
  &&
  \left\{
  \begin{array}{@{\hspace{0mm}}c@{\hspace{0mm}}}
    e^{1 1}_{1 1}
  \end{array}
  \right\},
  \qquad
  -
  p^2 \overline{Q}^2
  \left\{
  \begin{array}{@{\hspace{0mm}}c@{\hspace{0mm}}}
    e^{2 2}_{2 2}
    \\
    e^{3 3}_{3 3}
  \end{array}
  \right\},
  \qquad
  p^4
  \left\{
  \begin{array}{@{\hspace{0mm}}c@{\hspace{0mm}}}
    e^{4 4}_{4 4}
  \end{array}
  \right\},
  \\
  &&
  -
  p \overline{Q}
  (p \overline{Q} - \overline{p} Q)
  \left\{
  \begin{array}{@{\hspace{0mm}}c@{\hspace{0mm}}}
    e^{2 1}_{2 1}
  \end{array}
  \right\},
  \qquad
  p^3 \overline{Q}
  (p Q - \overline{p} \overline{Q})
  \left\{
  \begin{array}{@{\hspace{0mm}}c@{\hspace{0mm}}}
    e^{4 3}_{4 3}
  \end{array}
  \right\},
  \\
  &&
  p \overline{Q}
  \left\{
  \begin{array}{@{\hspace{0mm}}c@{\hspace{0mm}}}
    e^{1 2}_{2 1},
    e^{1 3}_{3 1}
    \\
    e^{2 1}_{1 2},
    e^{3 1}_{1 3}
  \end{array}
  \right\},
  \qquad
  p^3 \overline{Q}
  \left\{
  \begin{array}{@{\hspace{0mm}}c@{\hspace{0mm}}}
    e^{2 4}_{4 2},
    e^{3 4}_{4 3}
    \\
    e^{4 2}_{2 4},
    e^{4 3}_{3 4}
  \end{array}
  \right\},
  \qquad
  p^2 \overline{Q}^2
  \left\{
  \begin{array}{@{\hspace{0mm}}c@{\hspace{0mm}}}
    e^{1 4}_{4 1}
    \\
    e^{4 1}_{1 4}
  \end{array}
  \right\},
  \qquad
  -
  p^2
  \left\{
  \begin{array}{@{\hspace{0mm}}c@{\hspace{0mm}}}
    e^{2 3}_{3 2}
    \\
    e^{3 2}_{2 3}
  \end{array}
  \right\},
  \\
  &&
  -
  p^2 \overline{Q}
  {(p\overline{Q} - \overline{p}Q)}^{\frac{1}{2}}
  {(p Q - \overline{p} \overline{Q})}^{\frac{1}{2}}
  \left\{
  \begin{array}{@{\hspace{0mm}}c@{\hspace{0mm}}}
    e^{4 1}_{2 3}
    \\
    e^{2 3}_{4 1}
  \end{array}
  \right\}.
\end{eqnarray*}

$\check{R}^3$ has $17$ nonzero components:
\begin{eqnarray*}
  &&
  \left\{
  \begin{array}{@{\hspace{0mm}}c@{\hspace{0mm}}}
    e^{1 1}_{1 1}
  \end{array}
  \right\},
  \qquad
  -
  p^2 \overline{Q}^2
  \left\{
  \begin{array}{@{\hspace{0mm}}c@{\hspace{0mm}}}
    e^{2 2}_{2 2}
    \\
    e^{3 3}_{3 3}
  \end{array}
  \right\},
  \qquad
  p^4
  \left\{
  \begin{array}{@{\hspace{0mm}}c@{\hspace{0mm}}}
    e^{4 4}_{4 4}
  \end{array}
  \right\},
  \\
  & &
  p^2
  (Q^2 - \overline{Q}^2)
  \left\{
  \begin{array}{@{\hspace{0mm}}c@{\hspace{0mm}}}
    e^{3 2}_{3 2}
  \end{array}
  \right\},
  \\
  &&
  p \overline{Q}
  \left\{
  \begin{array}{@{\hspace{0mm}}c@{\hspace{0mm}}}
    e^{1 2}_{2 1},
    e^{1 3}_{3 1}
    \\
    e^{2 1}_{1 2},
    e^{3 1}_{1 3}
  \end{array}
  \right\},
  \qquad
  p^3 \overline{Q}
  \left\{
  \begin{array}{@{\hspace{0mm}}c@{\hspace{0mm}}}
    e^{2 4}_{4 2},
    e^{3 4}_{4 3}
    \\
    e^{4 2}_{2 4},
    e^{4 3}_{3 4}
  \end{array}
  \right\},
  \qquad
  p^2 \overline{Q}^2
  \left\{
  \begin{array}{@{\hspace{0mm}}c@{\hspace{0mm}}}
    e^{1 4}_{4 1}
    \\
    e^{4 1}_{1 4}
  \end{array}
  \right\},
  \qquad
  -
  p^2
  \left\{
  \begin{array}{@{\hspace{0mm}}c@{\hspace{0mm}}}
    e^{2 3}_{3 2}
    \\
    e^{3 2}_{2 3}
  \end{array}
  \right\}.
\end{eqnarray*}

$\check{R}^4$ has $16$ nonzero components:
\begin{eqnarray*}
  &&
  \left\{
  \begin{array}{@{\hspace{0mm}}c@{\hspace{0mm}}}
    e^{1 1}_{1 1}
  \end{array}
  \right\},
  \qquad
  -
  p^2 \overline{Q}^2
  \left\{
  \begin{array}{@{\hspace{0mm}}c@{\hspace{0mm}}}
    e^{2 2}_{2 2}
    \\
    e^{3 3}_{3 3}
  \end{array}
  \right\},
  \qquad
  p^4
  \left\{
  \begin{array}{@{\hspace{0mm}}c@{\hspace{0mm}}}
    e^{4 4}_{4 4}
  \end{array}
  \right\},
  \\
  &&
  p \overline{Q}
  \left\{
  \begin{array}{@{\hspace{0mm}}c@{\hspace{0mm}}}
    e^{1 2}_{2 1},
    e^{1 3}_{3 1}
    \\
    e^{2 1}_{1 2},
    e^{3 1}_{1 3}
  \end{array}
  \right\},
  \qquad
  p^3 \overline{Q}
  \left\{
  \begin{array}{@{\hspace{0mm}}c@{\hspace{0mm}}}
    e^{2 4}_{4 2},
    e^{3 4}_{4 3}
    \\
    e^{4 2}_{2 4},
    e^{4 3}_{3 4}
  \end{array}
  \right\},
  \qquad
  p^2 \overline{Q}^2
  \left\{
  \begin{array}{@{\hspace{0mm}}c@{\hspace{0mm}}}
    e^{1 4}_{4 1}
    \\
    e^{4 1}_{1 4}
  \end{array}
  \right\},
  \qquad
  -
  p^2
  \left\{
  \begin{array}{@{\hspace{0mm}}c@{\hspace{0mm}}}
    e^{2 3}_{3 2}
    \\
    e^{3 2}_{2 3}
  \end{array}
  \right\}.
\end{eqnarray*}

Each $\check{R}^i$ has a distinct set of eigenvalues, with some
overlap.  Immediately, the construction of $\check{\mathcal{R}}(u)$
reminds us that the three distinct (diagonal) components
$(\check{R}^i)^{j j}_{j j}$ must be (gauge-independent) eigenvalues.
In the spectral limit $u\to\infty$, these become $1$, $-q^{2 \alpha}$
(twice), and $q^{4 \alpha + 2}$ respectively.  No matter what the
gauge, we will have a minimum of these $3$ eigenvalues; in fact
$\check{R}^1$ contains only these $3$.  Overall, we expect a maximum of
$15$ distinct eigenvalues; these are described in Table \ref{tab:eigs}.

\begin{table}[ht]
  \centering
  \begin{tabular}{c|*{8}{l}|{c}}
    Case ($i$) & \multicolumn{8}{c|}{Eigenvalues~of~$\check{R}^{i}$} & \# \\
    \hline
      & & & & & & & & & \\[-3mm]
    1 & $1$ & $- q^{2 \alpha}$ & $q^{4\alpha+2}$ & & & & & & $3$ \\
      & $1$ & $- p^{2} \overline{Q}^2$ & $p^{4}$           & & & & & & \\[1mm]
    \hline
      & & & & & & & & & \\[-3mm]
    2 & $1$ & $-q^{2 \alpha}$ & $q^{4\alpha+2}$ &
        $\pm q^{\alpha}$ & $\pm q^{3 \alpha+1}$ & & & & $7$ \\
      & $1$ & $- p^{2} \overline{Q}^2$ & $p^{4}$           &
        $\pm p \overline{Q}$    & $\pm p^{3} \overline{Q}$ & & & & \\[1mm]
    \hline
      & & & & & & & & & \\[-3mm]
    3 & $1$ & $-q^{2 \alpha}$ & $q^{4\alpha+2}$ &
        $\pm q^{\alpha}$       & $\pm q^{3 \alpha+1}$ & $q^{2 \alpha}$ &
        $q^{2 \alpha + 2}$ & & $9$ \\
      & $1$ & $-p^{2} \overline{Q}^2$ & $p^{4}$           &
        $\pm p     \overline{Q}$    & $\pm p^{3} \overline{Q}$  & $p^{2} \overline{Q}^2$&
        $p^{2} Q^{2}        $ & & \\[1mm]
    \hline
      & & & & & & & & & \\[-3mm]
    4 & $1$ & $-q^{2 \alpha}$ & $q^{4\alpha+2}$ &
        $\pm q^{\alpha}$       & $\pm q^{3 \alpha+1}$ & $q^{2 \alpha}$ &
      & $\pm q^{2 \alpha + 1}$ & $10$ \\
      & $1$ & $- p^{2} \overline{Q}^2$ & $p^{4}$           &
        $\pm p     \overline{Q}$     & $\pm p^{3} \overline{Q}$& $p^{2} \overline{Q}^2$ &
      & $\pm p^{2}$ & \\
  \end{tabular}
  \caption{Eigenvalues of quantum R matrices for various gauges.}
  \label{tab:eigs}
\end{table}

%%%%%%%%%%%%%%%%%%%%%%%%%%%%%%%%%%%%%%%%%%%%%%%%%%%%%%%%%%%%%%%%%%%%%%%%

\section{Ambient isotopy link invariants}

A state model for evaluation of link invariants of ambient isotopy
based on a quantum R matrix $\check{R}$, is defined by two parameters:
a representation of the braid generator $\sigma$, and a representation
of the `left handle' $C$ \cite{DeWit:99a}. In our case, we seek a
scaling factor $\kappa$ such that $\sigma^{\pm 1} = \kappa^{\pm 1}
\check{R}^{\pm 1}$, and a $4\times 4$ matrix $C$ such that (Einstein
summation convention):
\begin{equation}
  C^d_c
  \cdot
  (\sigma^{\pm 1})^{c a}_{d b}
  =
  \delta^a_b,
  \label{eq:R1}
\end{equation}
is satisfied.

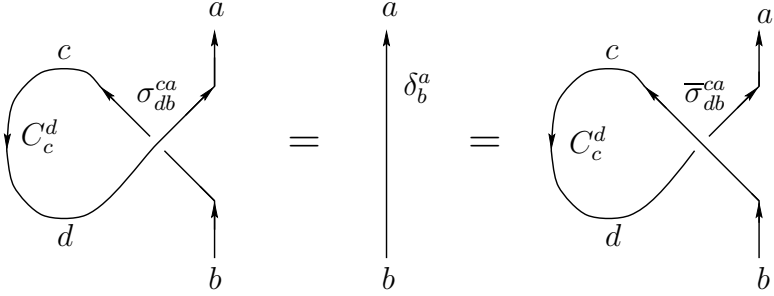
\begin{figure}[htbp]
  \begin{center}
    \input{graphics/Reidemeister1.pstex_t}
    \caption{
      The first Reidemeister move.
    }
    \label{fig:Reidemeister1}
  \end{center}
\end{figure}

This ensures that that the value of the link invariant over a
single loop of writhe is unity, viz that the first Reidemeister move
is satisfied.
This is depicted in Figure \ref{fig:Reidemeister1}.
In fact, this only determines $\kappa$ up to a
factor of $\pm 1$, but we shall ignore one case.

As $\check{R}$ necessarily satisfies the QYBE (\ref{eq:QYBEbraid}),
abstract tensors built from $\sigma$ are invariant under the second and
third Reidemeister moves, hence we may construct representations of
arbitrary braids from $\sigma$.  As all links may be represented by
braids combined with left handles, together these are sufficient
parameters.

Our state model then yields a polynomial invariant of oriented of
$(1,1)$ tangles.  Evaluation of such a model for arbitrary links is
described in our previous work \cite{DeWit:99a,DeWitKauffmanLinks:99a}.
Here, we use the same principles and \textsc{Mathematica} code -- only
the state model parameters are changed.

%%%%%%%%%%%%%%%%%%%%%%%%%%%%%%%%%%%%%%%%%%%%%%%%%%%%%%%%%%%%%%%%%%%%%%%%

\subsection*{Case 1}

Suitable $C$ that satisfy (\ref{eq:R1}) are:
\begin{eqnarray*}
  &&
  \overline{\kappa}
  \cdot
  \overline{p}^4
  Q^{2}
  A
  =
  \kappa
  \cdot
  \overline{Q}^2
  A,
  \qquad
  \mathrm{where}
  \qquad
  A =
  \mathrm{diag}
    \left\{
      + Q^{2},
      - Q^{2},
      - \overline{Q}^2,
      + \overline{Q}^2
    \right\},
\end{eqnarray*}
hence $\kappa=\overline{p}^{2}Q^{2}$ suffices to yield:
$ C = p^{2} A $.
The associated braid generator $\sigma$ is:
\begin{eqnarray*}
  &&
  \overline{p}^2 Q^2
  \left\{
  \begin{array}{@{\hspace{0mm}}c@{\hspace{0mm}}}
    e^{1 1}_{1 1}
  \end{array}
  \right\},
  \qquad
  -
  1
  \left\{
  \begin{array}{@{\hspace{0mm}}c@{\hspace{0mm}}}
    e^{2 2}_{2 2}
    \\
    e^{3 3}_{3 3}
  \end{array}
  \right\},
  \qquad
  p^2 Q^2
  \left\{
  \begin{array}{@{\hspace{0mm}}c@{\hspace{0mm}}}
    e^{4 4}_{4 4}
  \end{array}
  \right\},
  \\
  &&
  -
  \overline{p} Q
  (p \overline{Q} - \overline{p} Q)
  \left\{
  \begin{array}{@{\hspace{0mm}}c@{\hspace{0mm}}}
    e^{2 1}_{2 1}
    \\
    e^{3 1}_{3 1}
  \end{array}
  \right\},
  \qquad
  p Q
  (p Q - \overline{p} \overline{Q})
  \left\{
  \begin{array}{@{\hspace{0mm}}c@{\hspace{0mm}}}
    e^{4 2}_{4 2}
    \\
    e^{4 3}_{4 3}
  \end{array}
  \right\},
  \\
  &&
  Q^2
  (p \overline{Q} - \overline{p} Q)
  (p Q - \overline{p} \overline{Q})
  \left\{
  \begin{array}{@{\hspace{0mm}}c@{\hspace{0mm}}}
    e^{4 1}_{4 1}
  \end{array}
  \right\},
  \qquad
  Q^2
  ( Q^{2} - \overline{Q}^{2} )
  \left\{
  \begin{array}{@{\hspace{0mm}}c@{\hspace{0mm}}}
    e^{3 2}_{3 2}
  \end{array}
  \right\},
  \\
  &&
  \overline{p} Q
  \left\{
  \begin{array}{@{\hspace{0mm}}c@{\hspace{0mm}}}
    e^{1 2}_{2 1},
    e^{1 3}_{3 1}
    \\
    e^{2 1}_{1 2},
    e^{3 1}_{1 3}
  \end{array}
  \right\},
  \qquad
  p Q
  \left\{
  \begin{array}{@{\hspace{0mm}}c@{\hspace{0mm}}}
    e^{2 4}_{4 2},
    e^{3 4}_{4 3}
    \\
    e^{4 2}_{2 4},
    e^{4 3}_{3 4}
  \end{array}
  \right\},
  \qquad
  \left\{
  \begin{array}{@{\hspace{0mm}}c@{\hspace{0mm}}}
    e^{1 4}_{4 1}
    \\
    e^{4 1}_{1 4}
  \end{array}
  \right\},
  \qquad
  -
  Q^2
  \left\{
  \begin{array}{@{\hspace{0mm}}c@{\hspace{0mm}}}
    e^{2 3}_{3 2}
    \\
    e^{3 2}_{2 3}
  \end{array}
  \right\},
  \\
  &&
  Q^2
  {(p\overline{Q} - \overline{p}Q)}^{\frac{1}{2}}
  {(p Q - \overline{p} \overline{Q})}^{\frac{1}{2}}
  \left\{
  \begin{array}{@{\hspace{0mm}}c@{\hspace{0mm}}}
    -\overline{Q}
    \left\{
    \begin{array}{@{\hspace{0mm}}c@{\hspace{0mm}}}
      e^{4 1}_{2 3}
      \\
      e^{2 3}_{4 1}
    \end{array}
    \right\},
    Q
    \left\{
    \begin{array}{@{\hspace{0mm}}c@{\hspace{0mm}}}
      e^{3 2}_{4 1}
      \\
      e^{4 1}_{3 2}
    \end{array}
    \right\}
  \end{array}
  \right\}.
\end{eqnarray*}
These choices of course lead to the Links--Gould invariant $LG$ of
\cite{DeWit:99a,DeWitKauffmanLinks:99a}.

For Cases 2--4, there are problems in finding $\kappa$, as the process
yields no consistent choices of $\kappa$ without reducing the variables
$p$ and $Q$.  (Furthermore, in each case, relaxing the condition that
$C$ be diagonal immediately yields the conclusion that $C$ \emph{must}
be diagonal if it does exist.)

\subsection*{Case 2}

Appropriate $C$ are:
\begin{eqnarray*}
  &&
  \overline{\kappa}
  \cdot
  \overline{p}^3
  Q
  \cdot
  \mathrm{diag}
  \left\{
    + p Q,
    - p Q,
    - \overline{p} \overline{Q},
    + \overline{p} \overline{Q}
  \right\},
  \\
  &&
  \kappa
  \cdot
  p
  \overline{Q}
  \cdot
  \mathrm{diag}
  \left\{
    + \overline{p} Q,
    - \overline{p} Q,
    - p \overline{Q},
    + p \overline{Q}
  \right\}.
\end{eqnarray*}
The only solution to this system is found by setting $p=\pm 1$, whence
$\kappa = Q$. In this case, we have:
\begin{eqnarray*}
  C
  =
  \mathrm{diag}
  \left\{
    + Q,
    - Q,
    - \overline{Q},
    + \overline{Q}
  \right\}
\end{eqnarray*}
and the following one-variable braid generator $\sigma$ (note the
imaginary components):
\begin{eqnarray*}
  &&
  Q
  \left\{
  \begin{array}{@{\hspace{0mm}}c@{\hspace{0mm}}}
    e^{1 1}_{1 1}
  \end{array}
  \right\},
  \qquad
  -
  \overline{Q}
  \left\{
  \begin{array}{@{\hspace{0mm}}c@{\hspace{0mm}}}
    e^{2 2}_{2 2}
    \\
    e^{3 3}_{3 3}
  \end{array}
  \right\},
  \qquad
  Q
  \left\{
  \begin{array}{@{\hspace{0mm}}c@{\hspace{0mm}}}
    e^{4 4}_{4 4}
  \end{array}
  \right\},
  \\
  &&
  (Q - \overline{Q})
  \left\{
  \begin{array}{@{\hspace{0mm}}c@{\hspace{0mm}}}
    e^{2 1}_{2 1},
    e^{4 3}_{4 3}
  \end{array}
  \right\},
  \\
  &&
  \pm
  \left\{
  \begin{array}{@{\hspace{0mm}}c@{\hspace{0mm}}}
    e^{1 2}_{2 1},
    e^{1 3}_{3 1}
    \\
    e^{2 1}_{1 2},
    e^{3 1}_{1 3}
  \end{array}
  \right\},
  \qquad
  \pm
  \left\{
  \begin{array}{@{\hspace{0mm}}c@{\hspace{0mm}}}
    e^{2 4}_{4 2},
    e^{3 4}_{4 3}
    \\
    e^{4 2}_{2 4},
    e^{4 3}_{3 4}
  \end{array}
  \right\},
  \qquad
  \overline{Q}
  \left\{
  \begin{array}{@{\hspace{0mm}}c@{\hspace{0mm}}}
    e^{1 4}_{4 1}
    \\
    e^{4 1}_{1 4}
  \end{array}
  \right\},
  \qquad
  -
  Q
  \left\{
  \begin{array}{@{\hspace{0mm}}c@{\hspace{0mm}}}
    e^{2 3}_{3 2}
    \\
    e^{3 2}_{2 3}
  \end{array}
  \right\},
  \\
  &&
  -
  i
  (Q - \overline{Q})
  \left\{
  \begin{array}{@{\hspace{0mm}}c@{\hspace{0mm}}}
    e^{4 1}_{2 3}
    \\
    e^{2 3}_{4 1}
  \end{array}
  \right\}.
\end{eqnarray*}

After the scaling, the $7$ distinct eigenvalues of $\check{R}^2$
coalesce to $4$ of $\sigma$; they are
$
  \left\{
      q^{2},
    - \overline{q}^2
    \pm 1
  \right\}
  \equiv
  \left\{
      Q,
    - \overline{Q},
    \pm 1
  \right\}
$.
Inspection of the eigenvectors shows that the $\sigma$ is not
diagonizable, as it has only $14$ (out of $16$) distinct eigenvectors.
%%% (Compare this with $15$ in Case 1 in the limit $p=1$.)

The resulting invariant is nontrivial.  It turns out to be the
Alexander--Conway invariant in variable $q$. (This is an experimental
observation only, but has been verified for all prime knots of up to
$10$ crossings.)

%%%%%%%%%%%%%%%%%%%%%%%%%%%%%%%%%%%%%%%%%%%%%%%%%%%%%%%%%%%%%%%%%%%%%%%%

\subsection*{Case 3}

Appropriate $C$ are:
\begin{eqnarray*}
  &&
  \overline{\kappa}
  \cdot
  \overline{p}^2
  Q^2
  \cdot
  \mathrm{diag}
  \left\{
    + p^2 \overline{Q}^2,
    - Q^{4},
    - \overline{Q}^{4},
    + \overline{p}^{2} \overline{Q}^2
  \right\}
  \quad
  \textrm{and}
  \\
  &&
  \kappa
  \cdot
  p^2
  \overline{Q}^2
  \cdot
  \mathrm{diag}
  \left\{
    + \overline{p}^2 Q^2,
    - Q^4,
    - \overline{Q}^{4},
    + p^{2} Q^2
  \right\}.
\end{eqnarray*}
The only solutions to this system leave us with an integer invariant
(e.g. $p=Q=1$ and $\kappa=1$). In any case, setting $Q=1$ means that the
R matrix degenerates to being a special case of Case 4, so we ignore
these solutions.

%%%%%%%%%%%%%%%%%%%%%%%%%%%%%%%%%%%%%%%%%%%%%%%%%%%%%%%%%%%%%%%%%%%%%%%%

\subsection*{Case 4}

Appropriate $C$ are:
\begin{eqnarray*}
  & &
  \overline{\kappa}
  \cdot
  \overline{p}^{2}
  \cdot
  \mathrm{diag}
  \left\{
            p^{2},
   - Q^{2},
   - Q^{2},
   +        \overline{p}^{2}
  \right\}
  \quad
  \textrm{and}
  \\
  & &
  \kappa
  \cdot
  p^{2}
  \cdot
  \mathrm{diag}
  \left\{
            \overline{p}^{2},
   - \overline{Q}^2,
   - \overline{Q}^2,
   +        p^{2}
  \right\}.
\end{eqnarray*}

Again, the only solutions to this system leave us with an integer
invariant. For example, for $p=Q=\pm 1$ and $\kappa=1$, we obtain
$
  C
  =
  \mathrm{diag}
  \left\{
    +1, -1, -1, +1
  \right\}
$,
and the braid generator:
\begin{eqnarray*}
  &&
  \left\{
  \begin{array}{@{\hspace{0mm}}c@{\hspace{0mm}}}
    e^{1 1}_{1 1}
  \end{array}
  \right\},
  \qquad
  -
  \left\{
  \begin{array}{@{\hspace{0mm}}c@{\hspace{0mm}}}
    e^{2 2}_{2 2}
    \\
    e^{3 3}_{3 3}
  \end{array}
  \right\},
  \qquad
  \left\{
  \begin{array}{@{\hspace{0mm}}c@{\hspace{0mm}}}
    e^{4 4}_{4 4}
  \end{array}
  \right\},
  \\
  &&
  \pm
  \left\{
  \begin{array}{@{\hspace{0mm}}c@{\hspace{0mm}}}
    e^{1 2}_{2 1},
    e^{1 3}_{3 1}
    \\
    e^{2 1}_{1 2},
    e^{3 1}_{1 3}
  \end{array}
  \right\},
  \qquad
  \pm
  \left\{
  \begin{array}{@{\hspace{0mm}}c@{\hspace{0mm}}}
    e^{2 4}_{4 2},
    e^{3 4}_{4 3}
    \\
    e^{4 2}_{2 4},
    e^{4 3}_{3 4}
  \end{array}
  \right\},
  \qquad
  \left\{
  \begin{array}{@{\hspace{0mm}}c@{\hspace{0mm}}}
    e^{1 4}_{4 1}
    \\
    e^{4 1}_{1 4}
  \end{array}
  \right\},
  \qquad
  -
  \left\{
  \begin{array}{@{\hspace{0mm}}c@{\hspace{0mm}}}
    e^{2 3}_{3 2}
    \\
    e^{3 2}_{2 3}
  \end{array}
  \right\}.
\end{eqnarray*}
The $10$ distinct eigenvalues of $\check{R}^4$ coalesce to $2$ of
$\sigma$; they are just $ \left\{ \pm 1 \right\} $, and $\sigma$ is of
course diagonizable.

The resulting invariant is integer, and in fact trivial.  This has been
confirmed by the application of the `Matveev $\Delta$--$\nabla$ test'
\cite{Matveev:87}. This test applies a little theorem which shows that
if an invariant cannot distinguish the braids
$\sigma_1\overline{\sigma}_2\sigma_1$ and
$\sigma_2\overline{\sigma}_1\sigma_2$, then it cannot distinguish any
knot from the unknot, since repeated interchanges of these braids are
sufficient to untangle any knot.

In fact, applying this test to $\check{R}^4$ shows that, regardless of
$C$ or any special choice of the variables $p$ and $Q$, any associated
invariant will be trivial. (Applying the test to Cases 2 and 3 shows
that we \emph{may} have a nontrivial invariant.)

%%%%%%%%%%%%%%%%%%%%%%%%%%%%%%%%%%%%%%%%%%%%%%%%%%%%%%%%%%%%%%%%%%%%%%%%

\section{Regular isotopy link invariants}

For Cases 2 and 3, in general, we fail to build a link invariant of
ambient isotopy as we cannot satisfy (\ref{eq:R1}).  We now build
invariants associated with these cases that are only of regular
isotopy.

\subsection*{Case 2}

The choice $\kappa=\overline{p}^{2} Q$ and
$
  C
  =
  \overline{p}
  \cdot
  \mathrm{diag}
  \left\{
    +Q, -Q, -\overline{Q}, +\overline{Q}
  \right\}
$,
yields the symmetric results:
\begin{eqnarray*}
  (\mathrm{tr} \otimes I)
  \;
  (C \otimes I)
  \;
  \sigma^{\pm 1}
  =
  \mathrm{diag}
  \left\{
    \overline{p}, \,
    \overline{p}, \,
    p, \,
    p
  \right\}^{\pm 1}.
\end{eqnarray*}
Applying our link invariant evaluation engine to these parameters, we
obtain an open tangle invariant which is a diagonal $4\times 4$
matrix.  For a knot $K$ (presented as the closure $\hat{\beta}$ of a
braid $\beta$), we obtain the corresponding \emph{regularly isotopic}
$(1,1)$-tangle invariant of the following form:
\begin{eqnarray*}
  \mathrm{diag}
  \{
    \underbrace{
      \overline{p}^{w}
      \Delta_K (Q^2 \overline{p}^2)
    }_{\textrm{(twice)}},
    \;
    \underbrace{
      p^{w}
      \Delta_K (Q^2 p^2)
    }_{\textrm{(twice)}}
  \},
\end{eqnarray*}
where $w$ is the writhe of $\beta$, and $\Delta_K$ is the
Alexander--Conway invariant of $K$. Again, this result is an
experimental observation only, checked to be valid for all prime knots
of up to $10$ crossings. Setting $p=1$ recovers our previous
observation that the only possible invariant of ambient isotopy is
$\Delta_K$.

%%%%%%%%%%%%%%%%%%%%%%%%%%%%%%%%%%%%%%%%%%%%%%%%%%%%%%%%%%%%%%%%%%%%%%%%

\subsection*{Case 3}

The choices $\kappa = \overline{p}^{2}$ and
$
  C
  =
  \mathrm{diag}
  \left\{
    1, Q^{2}, \overline{Q}^2, 1
  \right\}
$,
yield the symmetric results:
\begin{eqnarray*}
  (\mathrm{tr} \otimes I)
  \;
  (C \otimes I)
  \;
  \sigma^{\pm 1}
  =
  \mathrm{diag}
  \left\{
    + \overline{p}^2, \;
    -        \overline{Q}^{4}, \;
    -        \overline{Q}^{4}, \;
    + p^{2}
  \right\}^{\pm 1}.
\end{eqnarray*}
Again, we obtain an open-tangle invariant which is a diagonal
$4\times 4$ matrix:
\begin{eqnarray*}
  \mathrm{diag}
  \{
    \overline{p}^{2 w},
    \;
    \underbrace{
      (-\overline{Q}^4)^w
      V_K (Q^4)
    }_{\textrm{twice}},
    \;
    p^{2 w}
  \},
\end{eqnarray*}
where $V_K$ is the Jones polynomial of the link $K$. Again, we
emphasise that this is an experimental observation, known to be valid
for all prime knots of up to $10$ crossings.

%%%%%%%%%%%%%%%%%%%%%%%%%%%%%%%%%%%%%%%%%%%%%%%%%%%%%%%%%%%%%%%%%%%%%%%%

\section{Conclusions}

Our calculations show that it is possible to access several invariants
from a single solution of the TYBE. Specifically, we have shown that it
is possible to recover both the Jones and the Alexander--Conway
invariants from the solution of the TYBE which was originally employed
to define the Links--Gould invariant $LG$
\cite{DeWit:99a,DeWitKauffmanLinks:99a}.

Repeating this process for other solutions may well uncover hitherto
unknown invariants!  Particularly, new solutions to the TYBE
recently reported in \cite{DeWit:99d,DeWit:99e} warrant
investigation in this context.

%%%%%%%%%%%%%%%%%%%%%%%%%%%%%%%%%%%%%%%%%%%%%%%%%%%%%%%%%%%%%%%%%%%%%%%%

\section*{Acknowledgements}

Jon Links is grateful to the Australian Research Council for financial
support and the cherry blossoms of ``Hanami'' in Kyoto, April 2000 for
inspiration.

David De Wit's research at Kyoto University is funded by a Postdoctoral
Fellowship for Foreign Researchers (\# P99703), provided by the Japan
Society for the Promotion of Science.
%%% $B$I$&$b(B $B$"$j$,$H$&(B $B$4$6$$$^$7$?!#(B
D\={o}mo arigat\={o} gozaimashita!

%%%%%%%%%%%%%%%%%%%%%%%%%%%%%%%%%%%%%%%%%%%%%%%%%%%%%%%%%%%%%%%%%%%%%%%%

\bibliographystyle{plain}
\bibliography{LinksDeWit2000}

\end{document}

%% file: graphics/Reidemeister1.pstex_t
\begin{picture}(0,0)%
\epsfig{file=graphics/Reidemeister1.pstex}%
\end{picture}%
\setlength{\unitlength}{0.00050000in}%
\begingroup\makeatletter\ifx\SetFigFont\undefined%
\gdef\SetFigFont#1#2#3#4#5{%
  \reset@font\fontsize{#1}{#2pt}%
  \fontfamily{#3}\fontseries{#4}\fontshape{#5}%
  \selectfont}%
\fi\endgroup%
\begin{picture}(7919,3045)(504,-2800)
\put(1051,-361){\makebox(0,0)[lb]{\smash{\SetFigFont{7}{8.4}{\familydefault}{\mddefault}{\updefault}{\large $c$}}}}
\put(676,-1261){\makebox(0,0)[lb]{\smash{\SetFigFont{7}{8.4}{\familydefault}{\mddefault}{\updefault}{\large $C^d_c$}}}}
\put(3451,-1336){\makebox(0,0)[lb]{\smash{\SetFigFont{7}{8.4}{\familydefault}{\mddefault}{\updefault}{\LARGE $=$}}}}
\put(5326,-1336){\makebox(0,0)[lb]{\smash{\SetFigFont{7}{8.4}{\familydefault}{\mddefault}{\updefault}{\LARGE $=$}}}}
\put(2626, 89){\makebox(0,0)[lb]{\smash{\SetFigFont{7}{8.4}{\familydefault}{\mddefault}{\updefault}{\large $a$}}}}
\put(4426, 89){\makebox(0,0)[lb]{\smash{\SetFigFont{7}{8.4}{\familydefault}{\mddefault}{\updefault}{\large $a$}}}}
\put(6751,-361){\makebox(0,0)[lb]{\smash{\SetFigFont{7}{8.4}{\familydefault}{\mddefault}{\updefault}{\large $c$}}}}
\put(2626,-2761){\makebox(0,0)[lb]{\smash{\SetFigFont{7}{8.4}{\familydefault}{\mddefault}{\updefault}{\large $b$}}}}
\put(4426,-2761){\makebox(0,0)[lb]{\smash{\SetFigFont{7}{8.4}{\familydefault}{\mddefault}{\updefault}{\large $b$}}}}
\put(8326,-2761){\makebox(0,0)[lb]{\smash{\SetFigFont{7}{8.4}{\familydefault}{\mddefault}{\updefault}{\large $b$}}}}
\put(6751,-2311){\makebox(0,0)[lb]{\smash{\SetFigFont{7}{8.4}{\familydefault}{\mddefault}{\updefault}{\large $d$}}}}
\put(1051,-2311){\makebox(0,0)[lb]{\smash{\SetFigFont{7}{8.4}{\familydefault}{\mddefault}{\updefault}{\large $d$}}}}
\put(6376,-1336){\makebox(0,0)[lb]{\smash{\SetFigFont{7}{8.4}{\familydefault}{\mddefault}{\updefault}{\large $C^d_c$}}}}
\put(4651,-736){\makebox(0,0)[lb]{\smash{\SetFigFont{7}{8.4}{\familydefault}{\mddefault}{\updefault}{\large $\delta^a_b$}}}}
\put(8326, 14){\makebox(0,0)[lb]{\smash{\SetFigFont{7}{8.4}{\familydefault}{\mddefault}{\updefault}{\large $a$}}}}
\put(7576,-811){\makebox(0,0)[lb]{\smash{\SetFigFont{7}{8.4}{\familydefault}{\mddefault}{\updefault}{\large $\overline{\sigma}^{c a}_{d b}$}}}}
\put(1876,-811){\makebox(0,0)[lb]{\smash{\SetFigFont{7}{8.4}{\familydefault}{\mddefault}{\updefault}{\large $\sigma^{c a}_{d b}$}}}}
\end{picture}